\documentclass[final]{siamltex}
\usepackage{booktabs}
\usepackage{multirow}
\usepackage{amsmath,amsfonts}
\usepackage{colortbl}
\usepackage{hyperref}
\usepackage{algorithmic, algorithm}

\newtheorem{remark}[theorem]{Remark}

\definecolor{gray}{rgb}{0.80,0.80,0.80}
\definecolor{lightgray}{rgb}{0.92,0.92,0.92}

\title{Differential qd algorithm with shifts for rank-structured matrices}

\author{Pavel Zhlobich\thanks{School of Mathematics, The University of Edinburgh, JCMB, The King's Buildings, Edinburgh, Scotland EH9 3JZ ({\tt P.Zhlobich@ed.ac.uk}).}}

\begin{document}

\maketitle

\begin{abstract}
Although QR iterations dominate in eigenvalue computations, there are several important cases when alternative LR-type algorithms may be preferable. In particular, in the symmetric tridiagonal case where differential qd algorithm with shifts (dqds) proposed by Fernando and Parlett enjoys often faster convergence while preserving high relative accuracy. In eigenvalue computations for rank-structured matrices QR algorithm is also a popular choice since, in the symmetric case, the rank structure is preserved. In the unsymmetric case, however, QR algorithm destroys the rank structure and, hence, LR-type algorithms come to play once again. In the current paper we adapt several variants of qd algorithms to quasiseparable matrices. Remarkably, one of them, when applied to Hessenberg matrices, becomes a direct generalization of dqds algorithm for tridiagonal matrices. Therefore, it can be applied to such important matrices as companion and confederate, and provides an alternative algorithm for finding roots of a polynomial represented in a basis of orthogonal polynomials. Results of preliminary numerical experiments are presented.
\end{abstract}

\begin{keywords} 
quasiseparable matrix, eigenvalue problem, dqds algorithm, companion matrix, comrade matrix
\end{keywords}

\begin{AMS}
15A18, 15A23, 65F15
\end{AMS}

\pagestyle{myheadings}
\thispagestyle{plain}
\markboth{P. ZHLOBICH}{DQDS ALGORITHM FOR RANK-STRUCTURED MATRICES}

\section{Introduction}\label{sec:intro}
Eigenvalue problem for rank-structured matrices such as semiseparable, quasiseparable, unitary-plus-rank-one and others has been an area of intense research in the last decade. This is due to the fact that the class of rank-structured matrices includes many other important classes, among which are banded matrices and their inverses, unitary Hessenberg matrices, companion and confederate matrices. Moreover, computations with such matrices are cheap and one step of any iterative algorithm usually takes only $\mathcal{O}(n)$ arithmetic operations, where $n$ is the size of the matrix.

There is one significant difference of rank-structured eigenvalue computations from the unstructured case. A prospective algorithm must preserve the low-rank structure of the initial matrix (maybe in some other form) to take advantage of fast linear algebra. This conservation property may not always be taken for granted and must be taken into account while developing new algorithms. 

As usual, there are two competitive approaches to eigenvalue ``hunting'': QR- and LR-type algorithms. Let us summarize the current state of the above mentioned approaches.
\begin{itemize}
\item[\textbf{QR:}] Development of QR-type eigenvalue solvers for rank-structured matrices was initially motivated by the matrix formulation of polynomial root-finding problem. Roots of a polynomial $P(x)=x^n+m_{n-1}x^{n-1}+\cdots+m_1x+m_0.$ are equal to the eigenvalues of its companion matrix:
\[
\begin{bmatrix}
-m_{n-1}  & -m_{n-2}  & \cdots & -m_1 & -m_0  \\
1     & 0     & \cdots & 0        & 0     \\
0     & 1     & \cdots & 0        & 0     \\
\cdot & \cdot & \cdot  & \cdot    & \cdot \\
0     & 0     & \cdots & 1        & 0     \\
\end{bmatrix}.
\]
This approach was successfully pursued in \cite{BDG04,BGP05} where certain low-rank preservation properties of QR-iterations for companion matrices were used. It was observed soon \cite{BEGG07,CGXZ08,BBEGG10} that more generally QR iterations preserve the structure of unitary-plus-rank-one matrices and companion matrix is only one particular representative of this class. Another direction of research aims at finding eigenvalues of symmetric quasiseparable matrices \cite{EGO05b,VVM05b}, the structure that is also preserved under QR iterations.

\item[\textbf{LR:}] One clear advantage of LR-type algorithms over the QR-type ones is that the former preserve quasiseparable structure of iterates even in the unsymmetric case. The disadvantage is the possible instability in floating point arithmetic. However, as it has been observed recently, even the use of orthogonal transformations for rank-structured matrices may lead to an unstable algorithm \cite{DOZ11}. There are two papers to be mentioned here: \cite{PVV08} gives a Cholesky LR algorithm for the symmetric positive definite quasiseparable matrices and \cite{BBD11} uses Neville representation of quasiseparable matrices to develop a qd-type method.
\end{itemize}

In the important case of symmetric positive definite tridiagonal matrices Fernando and Parlett \cite{FP94,P95} developed a root-free eigenvalue problem solver called \emph{differential qd algorithm with shifts}, or \emph{dqds} for short. The algorithm is fast and accurate and has become one of LAPACK's main eigenvalue routines. Although \cite{BBD11} attempted to transfer this algorithm to the quasiseparable case, the proposed algorithm is not dqds in the strict sense as it uses Neville-type eliminations. 

The main contribution of the current paper is a direct generalization of dqds algorithm of Fernando and Parlett to the quasiseparable matrix case. By achieving this goal we strictly follow the methodology of \cite{P95} and, as a by-product, derive several new eigenvalue algorithms applicable in different cases. We list below all the algorithms obtained in the paper.
\begin{itemize}
\item New LU decomposition algorithm for general quasiseparable matrices.
\item Stationary and Progressive qd algorithms with shifts for general quasiseparable matrices.
\item Differential qd algorithm with shifts (dqds) for Hessenberg quasiseparable matrices.
\item Specification of dqds algorithm for companion and comrade matrices.
\end{itemize}

In the symmetric positive-definite tridiagonal case both dqds algorithm of Fernando and Parlett and the variant of QR algorithm by Demmel and Kahan \cite{DK90} guarantee high relative accuracy of the computed eigenvalues (although dqds algorithm is not backward stable). The unsymmetric case is much more complex and, to the best of our knowledge, relative accuracy has not been proved for any of the existing algorithms even in the simple tridiagonal case. In this paper we are mainly interested in the unsymmetric eigenvalue problem and do not address the issue of relative accuracy. However, we can say definitely that the proposed dqds algorithm for Hessenberg quasiseparable matrices is \emph{not} backward stable as it is a straightforward generalization of tridiagonal dqds algorithm. Another issue is that the new algorithm uses LU factorization without pivoting at the initial step. We, therefore, assume that this factorization exists. Nevertheless, numerical experiments with many different matrices show that the new dqds algorithm often delivers more accurate result than its QR-based competitors.

The outline of the paper is as follows. In Section \ref{sec:qs} we recall the definition of quasiseparable matrices and also derive their LU factorization that will be our main tool later. A restricted version of this algorithm applicable to diagonal plus semiseparable matrices (a subclass of quasiseparable matrices) was derived in \cite{GKK85} and some formulae of it are similar to the ones in the inversion algorithm given in \cite{EG99a}. Still, new LU factorization algorithm, to the best of our knowledge, has never been published and may be useful to those who need a fast system solver for quasiseparable matrices. The new algorithm uses the idea of successive Schur complements and in this respect it is also different from the algorithm proposed in \cite[p.~171]{VVM08} that is based on the representation of the original matrix as a sum of lower and upper triangular matrices. Moreover, we provide an explicit pseudocode of the algorithm, while no pseudocode was provided in \cite{VVM08}. The complexity of the algorithm is $\mathcal{O}(N)$ and it is valid in the \emph{strongly regular} case (i.e. \emph{all its block leading principal minors are non-vanishing}). In the subsequent Section \ref{sec:qd} we present two versions of qd algorithm for general quasiseparable matrices: stationary and progressive. Section \ref{sec:dqds} is central in the paper and the new dqds algorithm is presented there. Section \ref{sec:GS} provides an alternative derivation of dqds algorithm via the generalized Gram--Schmidt process and also shed light on the meaning of some parameters arising in the algorithm. We next specialize the dqds algorithm for general Hessenberg quasiseparable matrices to more specific cases of companion and comrade matrices in Section \ref{sec:companion}. Results of preliminary numerical experiments are presented in the final Section \ref{sec:results}.

\section{Quasiseparable matrices}\label{sec:qs}
There are many important classes of structured matrices with the property of having low-rank blocks above and below the diagonal that one can meet in applications. Among the most well-known are semiseparable, quasiseparable, $H$-matrices, $H^2$-matrices and mosaic-skeleton matrices. In the current paper we are particularly interested in eigenvalue problem for quasiseparable matrices leaving possible extensions of the results to other classes of rank-structured matrices for future research. We will also consider the most general version of quasiseparable matrices, namely \emph{block quasiseparable matrices} as many of our results trivially generalize from the scalar to the block case. The usual definition of a block quasiseparable matrix is given below.
\begin{definition}[Rank definition of a block quasiseparable matrix]\label{def:qsrank}
Let $A$ be a block matrix of block sizes $\{n_k\}_{k=1}^n$ then it is called block $(r^l,r^u)$-quasiseparable if
\begin{eqnarray*}
&\max_{K}\,\rank A(K+1:N,1:K)\leq r^l,\quad\max_{K}\,\rank A(1:K,K+1:N)\leq r^u,\\
&K=\sum\limits_{i=1}^k n_i,\quad N=\sum\limits_{i=1}^n n_i,
\end{eqnarray*}
where $r^l$ ($r^u$) is called lower (upper) order of quasiseparability.
\end{definition}

In other words, quasiseparable matrices are those having low rank partitions in their upper and lower parts. In what follows we will call \emph{block $(r^l,r^u)$-quasiseparable} matrices simply \emph{quasiseparable} for shortness.

In order to exploit the quasiseparability of matrices in practice one must use a low-parametric representation of them. There are many alternative parametrizations of quasiseparable matrices all of which use $\mathcal{O}(N)$ parameters, where $N$ is the size of the matrix. Having such a parametrization at hand one can write most of the algorithms, e.g., inversion, LU, QR, matrix-vector multiplication in terms of these parameters and, therefore, the complexity of these algorithms is $\mathcal{O}(N)$. In the current paper we will use the so-called \emph{generator representation} (Definition \ref{def:qsgen} below) proposed by Eidelman and Gohberg \cite{EG99a}. For alternative parametrizations see \cite{XCGL09,DV07}.

\begin{definition}[Generator definition of a block quasiseparable matrix]\label{def:qsgen}
Let $A$ be a block matrix of block sizes $\{n_k\}_{k=1}^n$ then it is called block $(r^l,r^u)$-quasiseparable if it can be represented in the form
\begin{equation}\label{eq:qs_matrix}
\begin{bmatrix}
d_1 & g_1h_2 & g_1b_2h_3 & \cdots & g_1b_2\dots b_{n-1}h_n \\
p_2q_1 & d_2 & g_2h_3 & \cdots & g_2b_3\dots b_{n-1}h_n \\
p_3a_2q_1 & p_3q_2 & d_3 & \cdots & g_3b_4\dots b_{n-1}h_n \\
\vdots & \vdots & \vdots & \ddots & \vdots \\
p_na_{n-1}\dots a_2q_1 & p_na_{n-1}\dots a_3q_2 & p_na_{n-1}\dots
a_4q_3 & \cdots & d_n
\end{bmatrix},
\end{equation}
where parameters (called generators) $\{d_k$, $q_k$, $a_k$, $p_k$, $g_k$, $b_k$, $h_k\}$ are matrices of sizes as in the table below.

\begin{table}[H]
\begin{center}
\caption{Sizes of generators of a block quasiseparable matrix.}\label{tbl:uhgens}
\begin{tabular}[c]{l>{$}c<{$}>{$}c<{$}>{$}c<{$}>{$}c<{$}>{$}c<{$}>{$}c<{$}>{$}c<{$}}
& d_k & q_k & a_k & p_k & g_k & b_k & h_k\\
\midrule
\# of rows & n_k & r^l_k & r^l_k & n_{k} & n_k & r^u_{k-1} & r^u_{k-1} \\
\midrule
\# of cols & n_k & n_k & r^l_{k-1} & r^l_{k-1} & r^u_k & r^u_k & n^u_k \\
\bottomrule
\end{tabular}
\end{center}
\end{table}

Orders of quasiseparability $r^l$ and $r^u$ are maxima over the corresponding sizes of generators:
\[
r^l = \max_{1\leq k\leq n-1} r^l_k,\quad r^u = \max_{1\leq k\leq n-1} r^u_k.
\] 
\end{definition}

\begin{remark}
Generators are not unique, there are infinitely many ways to represent the same quasiseparable matrix with different generators. For the relation between different sets of generators see \cite{EG05}.
\end{remark}

\begin{remark}
Sizes $r^l_k$ and $r^u_k$ of generators are directly related to ranks of submatrices in the lower and upper parts correspondingly. Namely, if we let $K$ to be the block index: $K=\sum\limits_{i=1}^k n_i$, then
\begin{equation}\label{eq:subm_ranks}
\rank A(K+1:N,1:K) \leq r^l_k,\quad\rank A(1:K,K+1:N) \leq r^u_k.
\end{equation}
Moreover, for any quasiseparable matrix there exists a set of generators with sizes $r^l_k$ and $r^u_k$ that satisfy \eqref{eq:subm_ranks} with exact equalities (such generators are called \emph{minimal}). For their existence and construction see \cite{EG05}.
\end{remark}

We next derive LU factorization algorithm for a general block quasiseparable matrix in terms of the generators it is described by. First, let us note that quasiseparable structure of the original matrix implies the quasiseparable structure of $L$ and $U$ factors. The theorem below makes this statement precise and, in addition, relates generators of an original matrix to the generators of its factors.

\begin{theorem}\label{thm:LU}
Let $A$ be a strongly regular $N\times N$ block $(r^l,r^u)$-quasiseparable matrix given by generators $\{d_k,q_k,a_k,p_k,g_k,b_k,h_k\}$ as in \eqref{eq:qs_matrix}. Let $A=LU$ be its block LU decomposition of the same block sizes. Then
\begin{enumerate}
\item[(i)] Factors $L$ and $U$ are $(r^l,0)$-- and $(0,r^u)$-quasiseparable. Moreover, $r^l_k(L)=r^l_k(A)$ and $r^u_k(U)=r^u_k(A)$ for all $k=1,\ldots,n-1$.
\item[(ii)] $L$ and $U$ are parametrized by the generators
$\{I_k,\widetilde{q}_k,a_k,p_k,0,0,0\}$ and $\{\widetilde{d}_k$, $0$, $0$, $0$, $\widetilde{g}_k$, $b_k$, $h_k\}$, where $I_k$ are identity matrices of sizes $n_k\times n_k$ and new parameters $\widetilde{q}_k$, $\widetilde{d}_k$ and $\widetilde{g}_k$ can be computed using the following algorithm:
\begin{algorithm}[H]
\caption{Fast quasiseparable LU decomposition.}\label{alg:LU}
\begin{algorithmic}
\REQUIRE $d_k,q_k,a_k,p_k,g_k,b_k,h_k$
\STATE $\widetilde{d}_1=d_1,\quad \widetilde{q}_1=q_1\widetilde{d}_1^{-1},\quad \widetilde{g}_1=g_1,\quad f_1 = \widetilde{q}_1\widetilde{g}_1$
\FOR {$k=2$ \TO $n-1$}
\STATE $\widetilde{d}_k=d_k - p_kf_{k-1}h_k$
\STATE $\widetilde{q}_k=(q_k - a_kf_{k-1}h_k)\widetilde{d}_k^{-1}$
\STATE $\widetilde{g}_k=g_k - p_kf_{k-1}b_k$
\STATE $f_k = a_kf_{k-1}b_k+\widetilde{q}_k\widetilde{g}_k$
\ENDFOR
\STATE $\widetilde{d}_n=d_n - p_nf_{n-1}h_n.$
\ENSURE $\widetilde{d}_k,\widetilde{q}_k,\widetilde{g}_k$
\end{algorithmic}
\end{algorithm}
\end{enumerate}
\end{theorem}
\begin{proof}
Statement (i) of the theorem follows from statement (ii), so we only need to prove the latter part. 

Denote, as usual, by $K$ the block index: $K=\sum\limits_{i=1}^k n_i$ and note that quasiseparable representation \eqref{eq:qs_matrix} implies the following recurrence relation between the blocks of $A$:
\begin{equation}\label{eq:block_recurrences}
\begin{aligned}
&A=
\left[\begin{array}{c|c}
A(1:K,1:K) & G_kH_{k+1} \\
\hline
P_{k+1}Q_k & A(K+1:N,K+1:N)
\end{array}\right];\\
&Q_1 = q_1,\quad Q_k = [a_kQ_{k-1}\;\;q_k],\quad k = 2,\ldots n-1;\\
&P_n = p_n,\quad P_k = [p_k\;;\;P_{k+1}a_k],\quad k = n-1,\ldots 2;\\
&G_1 = g_1,\quad G_k = [G_{k-1}b_k\;;\;g_k],\quad k = 2,\ldots n-1;\\
&H_n = h_n,\quad H_k = [h_k\;\;b_kH_{k+1}],\quad k = n-1,\ldots 2.
\end{aligned}
\end{equation}

The proof will be constructed by induction. We will show that for each $k$
\begin{equation}\label{eq:blockLU}
\left[\begin{array}{c|c}
A_{11}^k & G_kH_{k+1} \\
\hline
P_{k+1}Q_k & \star
\end{array}\right]=
\left[\begin{array}{c|c}
L_{11}^k & 0 \\
\hline
P_{k+1}\widetilde{Q}_k & \star
\end{array}\right]
\cdot
\left[\begin{array}{c|c}
U_{11}^k & \widetilde{G}_kH_{k+1} \\
\hline
0 & \star
\end{array}\right].
\end{equation}

For $k=1$ we get from \eqref{eq:blockLU}:
\[
\begin{aligned}
&d_1 = \widetilde{d}_1,\\
&P_2Q_1=P_2\widetilde{Q}_1\widetilde{d}_1,\\
&G_1H_2 = \widetilde{G}_1H_2,
\end{aligned}\quad\Longleftarrow\quad
\begin{aligned}
&\widetilde{d}_1 = d_1,\\
&\widetilde{q}_1 = q_1\widetilde{d}_1^{-1},\\
&\widetilde{g}_1= g_1.
&\end{aligned}
\]

Let us introduce an auxiliary parameter $f_k=\widetilde{Q}_k\widetilde{G}_k$. It is easy to show by using \eqref{eq:block_recurrences} that $f_k$ satisfies the recurrence relation 
\begin{equation}\label{eq:f_k}
f_1 = \widetilde{q}_1\widetilde{g}_1,\quad f_k = a_kf_{k-1}b_k+\widetilde{q}_k\widetilde{g}_k.
\end{equation}

Assume that \eqref{eq:blockLU} holds for some fixed $k$, then it holds for $k+1$ if
\begin{eqnarray}
d_{k+1} = [p_{k+1}\widetilde{Q}_k \;\; 1]\cdot[\widetilde{G}_kh_{k+1} \;;\; \widetilde{d}_{k+1}],\label{eq1}\\
P_{k+2}q_{k+1} = P_{k+2}\widetilde{Q}_{k+1}\cdot[\widetilde{G}_kh_{k+1} \;;\; \widetilde{d}_{k+1}],\label{eq2}\\
g_{k+1}H_{k+2} = [p_{k+1}\widetilde{Q}_k \;\; 1]\cdot\widetilde{G}_{k+1}H_{k+2}.\label{eq3}
\end{eqnarray}

The first equality simplifies
\[
d_{k+1} = p_{k+1}\widetilde{Q}_k\widetilde{G}_kh_{k+1}+\widetilde{d}_{k+1} = p_{k+1}f_kh_{k+1}+\widetilde{d}_{k+1}.
\]
The second equality \eqref{eq2} holds if
\begin{multline*}
q_{k+1} = \widetilde{Q}_{k+1}\cdot[\widetilde{G}_kh_{k+1} \;;\; \widetilde{d}_{k+1}] = \\ = [a_{k+1}\widetilde{Q}_{k}\;\;\widetilde{q}_{k+1}]\cdot[\widetilde{G}_kh_{k+1} \;;\; \widetilde{d}_{k+1}] =
a_{k+1}f_kh_{k+1}+\widetilde{q}_{k+1}\widetilde{d}_{k+1}.
\end{multline*}
Finally, the last equality \eqref{eq3} is true if
\[
g_{k+1}=[p_{k+1}\widetilde{Q}_k \;\; 1]\cdot\widetilde{G}_{k+1} = 
[p_{k+1}\widetilde{Q}_k \;\; 1]\cdot[\widetilde{G}_{k}b_{k+1}\;;\;\widetilde{g}_{k+1}] = 
p_{k+1}f_kb_{k+1}+\widetilde{g}_{k+1}.
\]

Matrix $\widetilde{d}_{k+1}$ is invertible because, by our assumption, matrix $A$ is strongly regular. Hence, we conclude that \eqref{eq:blockLU} is true also for index $k+1$ if generators $\widetilde{q}_{k+1}$, $\widetilde{d}_{k+1}$ and $\widetilde{g}_{k+1}$ are those computed in Algorithm \ref{alg:LU}. The assertion of the theorem holds by induction.

\end{proof}

\section{qd algorithms for general quasiseparable matrices}\label{sec:qd}
In the previous section we have shown that factors in the quasiseparable LU factorization retain the low-rank structure of the original matrix. It turns out that even stronger result can be proved, namely that the low-rank structure is preserved under iterations of the LR algorithm. For completeness let us prove this simple theorem below.

\begin{theorem}\label{thm:LRinvariance}
Let $A$ be an $N\times N$ block $(r^l,r^u)$-quasiseparable matrix and $A - \sigma I$ be strongly regular. Define one step of shifted LR iterations by
\[
A - \sigma I = LU,\quad A' = UL.
\]
Then $A'$ is a strongly regular block $(r^l,r^u)$-quasiseparable matrix.
\end{theorem}
\begin{proof}
From Theorem \ref{thm:LU} we know that matrices $L$ and $U$ are $(r^l,0)$- and $(0,r^u)$-quasiseparable.
Let $\{n_k\}_{k=1}^n$ be sizes of blocks of $A$ and let $K=\sum\limits_{i=1}^k n_i$. Then for each $K$ let us write the product of $U$ and $L$ terms in the block form assuming that $A'_{11}=A'(1:K,1:K)$:
\[
\left[\begin{array}{c|c}
A'_{11} & A'_{12} \\
\hline
A'_{21} & A'_{22}
\end{array}\right]=
\left[\begin{array}{c|c}
U_{11} & U_{12} \\
\hline
0 & U_{22}
\end{array}\right]\cdot
\left[\begin{array}{c|c}
L_{11} & 0 \\
\hline
L_{21} & L_{22}
\end{array}\right].
\]
Hence,
\[
\begin{aligned}
&\rank(A'_{21})=\rank(L_{21}U_{11})=\rank(L_{21})=r^l_k,\\
&\rank(A'_{12})=\rank(U_{12}L_{22})=\rank(U_{12})=r^u_k,\\
\end{aligned}
\]
where $U_{11}$ is invertible by our assumption on the strong regularity of $A - \sigma I$.
\end{proof}

The assertion of Theorem \ref{thm:LRinvariance} lies in the heart of fast LR-type algorithms proposed in \cite{PVV08} and \cite{BBD11}. In fact Algorithm \ref{alg:LU} can be used to derive a new version of quasiseparable LR algorithm but we will not do it here as our main interest lies in deriving qd-type algorithms that are believed to be better in practice. Below we derive 3 new algorithms: \emph{stationary qd}, \emph{progressive qd} and \emph{Hessenberg dqds}. Our ultimate goal is the last algorithm and the first two can be seen as intermediate results.

\subsection{Stationary qd algorithm}
Triangular factors change in a complicated way under translation. Let $L$ and $U$ be quasiseparable factors as in Theorem \ref{thm:LU}, our task is to compute $\widehat{L}$ and $\widehat{U}$ so that
\[
A - \sigma I = LU - \sigma I = \widehat{L}\widehat{U}
\] 
for a given shift $\sigma$. Knowing that factors $L$ and $U$ on the left as well as $\widehat{L}$ and $\widehat{U}$ on the right are quasiseparable matrices we want to find formulae that define the direct mapping from the generators of the first pair to the generators of the second pair. Let factors $L$ and $U$ be given by generators $\{q_k$, $a_k$, $p_k\}$ and $\{d_k$, $g_k$, $b_k$, $h_k\}$, respectively. Simply inverting the formulae in Algorithm \ref{alg:LU} to get quasiseparable generators of $LU - \sigma I$ and $\widehat{L}\widehat{U}$ and equating them we obtain:
\[
\begin{gathered}
a_k = \widehat{a}_k,\quad b_k = \widehat{b}_k,\quad p_k = \widehat{p}_k,\quad h_k = \widehat{h}_k,\\
f_1 = 0,\quad f_{k+1} = a_kf_kb_k + q_kg_k,\quad\widehat{f}_1 = 0,\quad \widehat{f}_{k+1} = \widehat{a}_k\widehat{f}_k\widehat{b}_k + \widehat{q}_k\widehat{g}_k,\\
a_kf_kh_k + q_kd_k = \widehat{a}_k\widehat{f}_k\widehat{h}_k + \widehat{q}_k\widehat{d}_k,\quad p_kf_kb_k + g_k = \widehat{p}_k\widehat{f}_k\widehat{b}_k + \widehat{g}_k,\\
p_kf_kh_k + d_k - \sigma I_k = \widehat{p}_k\widehat{f}_k\widehat{h}_k + \widehat{d}_k.\\
\end{gathered}
\]
The above written formulae give rise to the algorithm of computing generators $\{\widehat{q}_k$, $\widehat{a}_k$, $\widehat{p}_k\}$ and $\{\widehat{d}_k$, $\widehat{g}_k$, $\widehat{b}_k$, $\widehat{h}_k\}$ of the new factors $\widehat{L}$ and $\widehat{U}$, respectively.

\begin{algorithm}[H]
\caption{Stationary $qd$ algorithm with shift (stqd($\sigma$)).}\label{alg:stat_qd}
\begin{algorithmic}
\REQUIRE $d_k$, $q_k$, $a_k$, $p_k$, $g_k$, $b_k$, $h_k$ and $\sigma$
\STATE $\widehat{d}_1=d_1-\sigma I_1,\quad\widehat{q}_1=q_1d_1\widehat{d}_1^{-1},\quad\widehat{g}_1=g_1$
\FOR {$k=2$ \TO $n-1$}
\STATE $\widehat{t}_{k} = a_{k-1}\widehat{t}_{k-1}b_{k-1}+q_{k-1}g_{k-1} - \widehat{q}_{k-1}\widehat{g}_{k-1}$
\STATE $\widehat{d}_k=p_k\widehat{t}_{k}h_k+d_k-\sigma I_k$
\STATE $\widehat{q}_k=\left(a_k\widehat{t}_kh_k + q_kd_k\right)\widehat{d}_k^{-1},\quad \widehat{g}_k=p_k\widehat{t}_kb_k + g_k$
\STATE $\widehat{a}_k=a_k,\quad\widehat{p}_k=p_k,\quad\widehat{b}_k=b_k,\quad\widehat{h}_k=h_k$
\ENDFOR
\STATE $\widehat{t}_n = a_{n-1}\widehat{t}_{n-1}b_{n-1}+q_{n-1}g_{n-1} - \widehat{q}_{n-1}\widehat{g}_{n-1}$
\STATE $\widehat{d}_n=p_n\widehat{t}_{n}h_n+d_n-\sigma I_n,\quad\widehat{p}_n=p_n,\quad\widehat{h}_n=h_n$
\ENSURE $\widehat{d}_k$, $\widehat{q}_k$, $\widehat{a}_k$, $\widehat{p}_k$, $\widehat{g}_k$, $\widehat{b}_k$, $\widehat{h}_k$
\end{algorithmic}
\end{algorithm}

Let us show that Algorithm \ref{alg:stat_qd} is the direct generalization of stationary $qd$ algorithm for tridiagonal matrices. $L$ and $U$ factors of a normalized\footnote{all entries in positions $(i,i+1)$ are ones} tridiagonal matrix with subdiagonal entries $l_k$ and diagonal entries $u_k$ have very special sets of quasiseparable generators, namely $\{1$, $l_k$, $0$, $1$, $0$, $0$, $0\}$ and $\{u_k$, $0$, $0$, $0$, $1$, $0$, $1\}$. For these special generators formulae in Algorithm \ref{alg:stat_qd} reduce to
\begin{equation}\label{eq:stat_qd_tridiag}
\begin{aligned}
&\widehat{u}_k \equiv \widehat{d}_k = p_k\widehat{t}_{k}h_k + d_k - \sigma = q_{k-1} + d_k - \sigma - \widehat{q}_{k-1} \equiv l_{k-1} + u_{k} - \sigma - \widehat{l}_{k-1},\\
&\widehat{l}_k \equiv \widehat{p}_{k+1}\widehat{q}_k = \widehat{p}_{k+1}\left(a_k\widehat{t}_kh_k + q_kd_k\right)\widehat{d}_k^{-1} = (p_{k+1}q_k)d_k\widehat{d}_k^{-1} \equiv l_ku_k/\widehat{u}_k.
\end{aligned}
\end{equation}
Equation \eqref{eq:stat_qd_tridiag} for $\widehat{l}_k$ and $\widehat{u}_k$ in terms of $l_k$ and $u_k$ is exactly the one defining stationary qd algorithm with shift for tridiagonal matrices (see \cite[page 465]{P95}).

\subsection{Progressive qd algorithm}

In deriving progressive qd algorithm with shift we seek for the triangular factorization of $A'-\sigma I$:
\begin{equation}\label{eq:LRiter}
A'-\sigma I = UL-\sigma I = \widehat{L}\widehat{U}
\end{equation}
for a suitable shift $\sigma$. Theorem \ref{thm:LRinvariance} tells us that $UL-\sigma I$ is a quasiseparable matrix of the same order as $L$ and $U$ factors and that suitable quasiseparable generators parametrization of $\widehat{L}$ and $\widehat{U}$ exists. However, we cannot obtain generators of $UL-\sigma I$ by simply inverting formulae in Algorithm \ref{alg:LU} as was done in the derivation of stationary qd algorithm. For this purpose we need to exploit the quasiseparable matrices multiplication Algorithm 4.3 from \cite{EG05} that defines generators of the product of two quasiseparable matrices.

Let $\{I_k$, $q_k$, $a_k$, $p_k$, $0$, $0$, $0\}$ and $\{d_k$, $0$, $0$, $0$, $g_k$, $b_k$, $h_k\}$ be generators of factors $L$ and $U$ in \eqref{eq:LRiter}, then generators  $\{d'_k$, $q'_k$, $a'_k$, $p'_k$, $g'_k$, $b'_k$, $h'_k\}$ of $A'-\sigma I = UL-\sigma I$ are as follows:
\begin{equation}\label{eq:ULprod}
\begin{aligned}
&q'_k = q_k,\quad g'_k = g_k,\quad a'_k = a_k,\quad b'_k = b_k,\\
&f_n = 0,\quad f_{k-1} = b_kf_ka_k + h_kp_k,\\
&p'_k = d_kp_k + g_kf_ka_k,\quad h'_k = h_k + b_kf_kq_k,\\
&d'_k = d_k + g_kf_kq_k - \sigma I_k.
\end{aligned}
\end{equation}

Equating generators in \eqref{eq:ULprod} to the generators of $\widehat{L}\widehat{U}$ we obtain progressive qd algorithm with shift described below.

\begin{algorithm}[H]
\caption{Progressive $qd$ algorithm with shift (qds($\sigma$)).}\label{alg:prog_qd}
\begin{algorithmic}
\REQUIRE $d_k$, $q_k$, $a_k$, $p_k$, $g_k$, $b_k$, $h_k$ and $\sigma$
\STATE $\widehat{p}_n=p_n,\quad \widehat{h}_n=h_n,\quad f_n=0$
\FOR {$k=n-1$ \TO $2$}
\STATE $f_k = b_{k+1}f_{k+1}a_{k+1}+h_{k+1}p_{k+1}$
\STATE $\widehat{p}_k = d_kp_k+g_kf_ka_k,\quad \widehat{h}_k = h_k+b_kf_kq_k$
\STATE $\widehat{a}_k = a_k,\quad\widehat{b}_k = b_k$
\ENDFOR
\STATE $f_1 = b_{2}f_{2}a_{2}+h_{2}p_{2},\quad \widehat{d}_1 = d_1+g_1f_1q_1-\sigma I_1$
\STATE $\widehat{q}_1 = q_1\widehat{d}_1^{-1},\quad \widehat{g}_1=g_1,\quad \widehat{f}_1 = 0$
\FOR {$k=2$ \TO $n-1$}
\STATE $\widehat{f}_{k} = \widehat{a}_{k-1}\widehat{f}_{k-1}\widehat{b}_{k-1}+\widehat{q}_{k-1}\widehat{g}_{k-1}$
\STATE $\widehat{q}_k = (q_k-\widehat{a}_k\widehat{f}_k\widehat{h}_k)\widehat{d}_k^{-1},\quad \widehat{g}_k = g_k-\widehat{p}_k\widehat{f}_k\widehat{b}_k$
\STATE $\widehat{d}_k = d_k+g_kf_kq_k-\sigma I_k-\widehat{p}_k\widehat{f}_k\widehat{h}_k$
\ENDFOR
\STATE $\widehat{f}_{n} = \widehat{a}_{n-1}\widehat{f}_{n-1}\widehat{b}_{n-1}+\widehat{q}_{n-1}\widehat{g}_{n-1},\quad \widehat{d}_n = d_n-\sigma I_n-\widehat{p}_n\widehat{f}_n\widehat{h}_n$
\ENSURE $\widehat{d}_k$, $\widehat{q}_k$, $\widehat{a}_k$, $\widehat{p}_k$, $\widehat{g}_k$, $\widehat{b}_k$, $\widehat{h}_k$
\end{algorithmic}
\end{algorithm}

In the case of normalized tridiagonal matrices, $L$ and $U$ factors have generators $\{1$, $l_k$, $0$, $1$, $0$, $0$, $0\}$ and $\{u_k$, $0$, $0$, $0$, $1$, $0$, $1\}$, and Algorithm \ref{alg:prog_qd} simplifies to
\begin{equation}\label{eq:prog_qd_tridiag}
\begin{aligned}
&\widehat{l}_k \equiv \widehat{p}_{k+1}\widehat{q}_k = (d_{k+1}p_{k+1})(q_k/\widehat{d}_k) \equiv l_ku_{k+1}/\widehat{u}_k,\\
&\widehat{u}_{k} \equiv \widehat{d}_{k} = d_{k}+p_{k+1}q_{k}-\sigma-\widehat{p}_{k}\widehat{q}_{k-1} \equiv u_{k}+l_{k}-\sigma-\widehat{l}_{k-1},
\end{aligned}
\end{equation}
where by $l_k/\widehat{l}_k$ and $u_k/\widehat{u}_k$ we denoted subdiagonal entries of $L/\widehat{L}$ factor and diagonal entries of $U/\widehat{U}$ factor correspondingly. One can easily spot in \eqref{eq:prog_qd_tridiag} the progressive $qd$ algorithm with shift for tridiagonal matrices \cite[page 467]{P95}.

\section{Differential qd algorithm in the Hessenberg quasiseparable case}\label{sec:dqds}
Fernando and Parlett \cite{FP94,P95} proposed a modified version of progressive qd algorithm called \emph{differential qd} algorithm. The proposed algorithm has received major attention due to its accuracy and speed, and is implemented as DLASQ in LAPACK. It turns out that differential version of progressive qd algorithm does not exist in the general case (the reason for this will be given in Section \ref{sec:GS}). However, we found that it exists for Hessenberg quasiseparable matrices that is an important case on its own as matrices related to polynomials, e.g. companion, confederate, fellow, are usually Hessenberg.

Consider a general Hessenberg quasiseparable matrix:
\begin{equation}\label{eq:Hqs_matrix}
\begin{bmatrix}
d_1 & g_1h_2 & g_1b_2h_3 & \cdots & g_1b_2\dots b_{n-1}h_n \\
s_1 & d_2 & g_2h_3 & \cdots & g_2b_3\dots b_{n-1}h_n \\
0 & s_2 & d_3 & \cdots & g_3b_4\dots b_{n-1}h_n \\
\vdots & \ddots & \ddots & \ddots & \vdots \\
0 & \cdots & 0 & s_{n-1} & d_n
\end{bmatrix}.
\end{equation}
Because of the Hessenberg structure, quasiseparable generators of the lower part are special, namely $q_k = s_k$, $a_k = 0$ and $p_k = 1$. Let us note that from now on we restrict matrices to the scalar case, i.e. generators such as $d_k$ are all scalars.

We next obtain differential version of Algorithm \ref{alg:prog_qd}. Our initial derivation may seem cumbersome and very technical but a much cleaner and justified derivation will be given later in Section \ref{sec:GS}.

First, observe that due to the special Hessenberg generators many formulae in Algorithm \ref{alg:prog_qd} can be simplified. In particular $f_k=h_{k+1}$ and $\widehat{f}_k=\widehat{q}_{k-1}\widehat{g}_{k-1}$. We next change the way $\widehat{d}_k$ are computed:
\[
\begin{aligned}
\widehat{d}_k &= d_k + s_kg_kh_{k+1} - \widehat{s}_{k-1}\widehat{g}_{k-1}\widehat{h}_k - \sigma\\
			  &= d_k + s_k(\widehat{g}_k +  \widehat{s}_{k-1}\widehat{g}_{k-1}b_k)h_{k+1} - \widehat{s}_{k-1}\widehat{g}_{k-1}\widehat{h}_k - \sigma\\
			  &= d_k + s_k\widehat{g}_kh_{k+1} - \widehat{s}_{k-1}\widehat{g}_{k-1}(\widehat{h}_k-s_kb_kh_{k+1}) - \sigma\\
			  &= d_k + s_k\widehat{g}_kh_{k+1} - \widehat{s}_{k-1}\widehat{g}_{k-1}h_k - \sigma.
\end{aligned}
\]
Next, let us define an auxiliary variable:
\begin{equation}\label{eq:aux}
t_k \equiv d_k - \widehat{s}_{k-1}\widehat{g}_{k-1}h_k - \sigma\;(=\widehat{d}_k-s_k\widehat{g}_kh_{k+1}).
\end{equation}
Observe that
\begin{multline}\label{eq:aux_rec}
t_{k+1} = d_{k+1}-(s_kd_{k+1}/\widehat{d}_k)\widehat{g}_{k}h_{k+1} - \sigma =\\= d_{k+1}/\widehat{d}_k(\widehat{d}_k-s_k\widehat{g}_{k}h_{k+1}) - \sigma = d_{k+1}t_k/\widehat{d}_k - \sigma.
\end{multline}
Using identities \eqref{eq:aux} and \eqref{eq:aux_rec} we can now derive dqds algorithm for Hessenberg quasiseparable matrices.

\begin{algorithm}[H]
\caption{Differential qd algorithm with shift (dqds($\sigma$)), Hessenberg case.}\label{alg:dqds_Hessenberg}
\begin{algorithmic}
\REQUIRE $s_k$, $d_k$, $g_k$, $b_k$, $h_k$ and $\sigma$
\STATE $t_1 = d_1 - \sigma,\quad \widehat{g}_1 = g_1,\quad \widehat{d}_1 = t_1 + s_1\widehat{g}_1h_2,$ 
\STATE $\widehat{s}_1=s_1d_2/\widehat{d}_1\quad t_2 = t_1d_2/\widehat{d}_1 - \sigma$
\FOR {$k=2$ \TO $n-1$}
\STATE $\widehat{h}_k = h_k+s_kb_kh_{k+1}$
\STATE $\widehat{g}_k = g_k - \widehat{s}_{k-1}\widehat{g}_{k-1}b_k$
\STATE $\widehat{d}_k = t_k + s_k\widehat{g}_kh_{k+1}$
\STATE $\widehat{s}_k = s_kd_{k+1}/\widehat{d}_k$
\STATE $t_{k+1} = t_{k}d_{k+1}/\widehat{d}_k - \sigma$
\ENDFOR
\STATE $\widehat{d}_n = t_n,\quad \widehat{h}_n = h_n$
\ENSURE $\widehat{s}_k$, $\widehat{d}_k$, $\widehat{g}_k$, $\widehat{h}_k$
\end{algorithmic}
\end{algorithm}

Let us show that Algorithm \ref{alg:dqds_Hessenberg} trivially reduces to the well-known dqds algorithm in the normalized tridiagonal case. In this case we have $b_k=0$ and $g_kh_{k+1}=1$ for all $k$ and, hence,
\[
\begin{aligned}
&\widehat{d}_k = t_k + s_k,\\
&\widehat{s}_k = s_kd_{k+1}/\widehat{d}_k,\\
&t_{k+1} = t_{k}(d_{k+1}/\widehat{d}_{k}) - \sigma.
\end{aligned}
\]
These formulae constitute dqds($\sigma$) algorithm for normalized tridiagonal matrices, see \cite[page 468]{P95}.

In the tridiagonal case dqds algorithm can be run from bottom to top of $L$ and $U$ factors. Let us show how Algorithm \ref{alg:dqds_Hessenberg} in the Hessenberg case can be also transformed to work backward. Let $A=LU$ be Hessenberg quasiseparable matrix with generators $\{s_k$, $d_k$, $g_k$, $b_k$, $h_k\}$ of $L$ and $U$ factors. Roots of characteristic polynomial $\det(x I-A)$ equal to the roots of $\det(x I-JAJ)$, where $J$ is antidiagonal matrix. Transform $JAJ$ is simply backward reordering of rows and columns of $A$. Applying this transform to the LU factorization we get $JAJ=(JLJ)(JUJ)$. It is trivial to show that matrices $(JLJ)$ and $(JUJ)$ retain all the properties of $L$ and $U$ and their quasiseparable generators are $\{s_{n-k}$, $d_{n-k+1}$, $h_{n-k+1}^T$, $b_{n-k+1}^T$, $g_{n-k+1}^T\}$. Hence, running Algorithm \ref{alg:dqds_Hessenberg} in a forward way on $JLJ$ and $JUJ$ is same as running it backward on $L$ and $U$.

\section{Generalized Gram--Schmidt process and quasiseparable dqds algorithm}\label{sec:GS}
Parlett \cite{P95} has shown that dqd algorithm for tridiagonal matrices can be interpreted as the generalized Gram--Schmidt orthogonalization process applied to matrices $L$ and $U$ from $A' = UL$. In this section we follow the same approach and derive Algorithm \ref{alg:dqds_Hessenberg} directly from the Gram--Schmidt and without referencing to Algorithms \ref{alg:prog_qd}. Our analysis uses little of quasiseparable matrices theory and is very accessible.

Let $F=\begin{bmatrix}f_1&f_2&\cdots&f_k\end{bmatrix}$ and $G=\begin{bmatrix}g_1&g_2&\cdots&g_k\end{bmatrix}$ be a pair of vector sets. Next theorem relates biorthogonalization of $F$ and $G$ (generalized Gram--Schmidt) to the LU factorization of $G^{*}F$ (here and thereafter $A^{*}$ denotes the conjugate transpose of $A$).

\begin{theorem}[Theorem 1, \cite{P95}]\label{thm:GS}
Let $F$ and $G$ be complex $n\times k$ matrices, $n\geq k$, such that $G^{*}F$ permits triangular factorization:
\[
G^{*}F=\widetilde{L}\widetilde{D}\widetilde{R},
\]
where $\widetilde{L}$ and $\widetilde{R}$ are unit triangular (left and right), respectively, and $\widetilde{D}$ is diagonal. Then there exist unique $n\times k$ matrices $\widetilde{Q}$ and $\widetilde{P}$ such that
\[
F=\widetilde{Q}\widetilde{R},\quad G=\widetilde{P}\widetilde{L}^{*},\quad \widetilde{P}^{*}\widetilde{Q}=\widetilde{D}.
\]
\end{theorem}

\begin{remark}
In practice, when $n=k$ and $\widetilde{D}$ is invertible one can omit $\widetilde{Q}$ and write $F = (\widetilde{P}^*)^{-1}(\widetilde{D}\widetilde{R})\equiv(\widetilde{P}^*)^{-1}\widetilde{U}$, $G = \widetilde{P}\widetilde{L}^*$ and still call it Gram--Schmidt factorization. The important feature is the uniqueness of $\widetilde{Q}$ and $\widetilde{P}$. The columns of $\widetilde{Q}$ and rows of $\widetilde{P}^{*}$ form a pair of dual bases for the space of $n$-vectors (columns) and its dual (the row $n$-vectors). There is no notion of orthogonality or inner product here; $p_i^{*}q_j=0$ simply says that $p_i^{*}$ annihilates $q_j$, $i\neq j$.
\end{remark}

To derive the unshifted version of dqds Algorithm \ref{alg:dqds_Hessenberg} apply Gram--Schmidt to the columns of $L$ and $U^*$, in the natural order, to obtain $\widehat{L}$ and $\widehat{U}$, then according to Theorem \ref{thm:GS}:
\[
UL = \widehat{L}\widehat{U}.
\]

Since matrix $L$ is unit lower bidiagonal, the matrix $\widehat{P}$ ($\widetilde{P}^*$ in Theorem \ref{thm:GS}) such that $L = \widehat{P}^{-1}\widehat{U}$ and $U = \widehat{L}\widehat{P}$ is upper Hessenberg. Fortunately, as will be shown soon, matrix $\widehat{P}$ has a simple $\mathcal{O}(n)$ representation versus dense $n(n-1)/2$ and makes the derivation of linear complexity algorithm possible. Matrices $L$ and $U$ will be transformed to $\widehat{U}$ and $\widehat{L}$ correspondingly by a sequence of simple transformations.

At the start of Gram--Schmidt factorization $L$ and $U$ factors are
\[
L=
\begin{bmatrix}
1 \\
s_1 & 1 \\
& s_2 & 1 \\
& & \ddots & \ddots 
\end{bmatrix},\quad
U=
\begin{bmatrix}
d_1 & g_1h_2 & g_1b_2h_3 & \cdots & g_1b_2\ldots b_{n-1}h_n \\
& d_2 & g_2h_3 & \cdots & g_2b_3\ldots b_{n-1}h_n \\
& & d_3 & \cdots & g_3b_4\ldots b_{n-1}h_n \\
& & & \ddots & \vdots
\end{bmatrix}.
\]

Let $\widehat{P} = P_{n-1}\ldots P_2P_1$. Choose the first transformation $P_1$ as
\[
P_1=
\begin{bmatrix}
d_1 & g_1h_2 & g_1b_2h_3 & \cdots & g_1b_2\ldots b_{n-1}h_n \\
-s_1 & 1 & 0 & \cdots & 0 \\
& & 1 & \ddots & \vdots \\
& & & \ddots & 0 \\
& & & & 1 \\
\end{bmatrix},
\]
then
\[
P_1^{-1}=
\begin{bmatrix}
1/\widehat{d}_1 & -g_1h_2/\widehat{d}_1 & -1/\widehat{d}_1\cdot g_1b_2h_3 & \cdots & -1/\widehat{d}_1\cdot g_1b_2\ldots b_{n-1}h_n \\
s_1/\widehat{d}_1 & d_1/\widehat{d}_1 & -s_1/\widehat{d}_1\cdot g_1b_2h_3 & \cdots & -s_1/\widehat{d}_1\cdot g_1b_2\ldots b_{n-1}h_n \\
& & 1 & \ddots & \vdots \\
& & & \ddots & 0 \\
& & & & 1 \\
\end{bmatrix},
\]
where $\widehat{d}_1=d_1+s_1g_1h_2$. The careful reader may check that matrices $P_1$ and $P_1^{-1}$ are such that
\[
P_1L=
\left[\begin{array}{c|cccl}
\widehat{d}_1 & g_1\widehat{h}_2 & g_1b_2\widehat{h}_3 & \cdots & g_1b_2\ldots b_{n-1}\widehat{h}_n \\
\hline
& 1 \\
& s_2 & 1 \\
& & s_3 & 1 \\
& & & \ddots & ~\ddots 
\end{array}\right]
\]
and
\[
UP_1^{-1}=
\left[\begin{array}{c|cccc}
1 & 0 & 0 & \cdots & 0 \\
\hline
\widehat{s}_1 & t_2 & \widehat{g}_2h_3 & \cdots & \widehat{g}_2b_3\ldots b_{n-1}h_n \\
& & d_3 & \cdots & g_3b_4\ldots b_{n-1}h_n \\
& & & \ddots & \vdots
\end{array}\right],
\]
where 
\begin{equation}\label{eq4}
\begin{aligned}
\widehat{d}_1 &=d_1+s_1g_1h_2,\quad \widehat{h}_k = h_k+s_kb_kh_{k+1},\\
\widehat{s}_1 &= s_1d_{2}/\widehat{d}_1,\quad \widehat{g}_2 = g_2 - \widehat{s}_{1}\widehat{g}_{1}b_2,\quad t_2 = d_{1}d_{2}/\widehat{d}_{1}.\\
\end{aligned}
\end{equation}

Note that $(2:n,2:n)$ submatrices of $P_1L$ and $UP_1^{-1}$ retained the form of the initial factors $L$ and $U$. Therefore, we can continue the orthogonalization process by using matrices $P_k$, shifted analogs of $P_1$. Formulae \eqref{eq4} as well as the ones we will obtain recursively repeat formulae of Algorithm \ref{alg:dqds_Hessenberg} in the case $\sigma = 0$. So, we have derived the dqd algorithm without referencing qd. The reason for the non-existence of the dqds algorithm in the general non-Hessenberg case is that if $L$ factor is dense, then matrices $P_k$ must also be dense and no recursive nested $\mathcal{O}(n)$ update of factor $U$ is possible.

One remarkable corollary of the Gram--Schmidt derivation of dqds algorithm is that it reveals the meaning of auxiliary quantities $t_k$.
\begin{theorem}\label{thm:t_k_noshift}
Let $A$ be a strongly regular Hessenberg quasiseparable matrix of an arbitrary order. Then quantities $t_k$ generated by Algorithm \ref{alg:dqds_Hessenberg} with zero shift are such that
\[
\frac{1}{t_k}=[A^{-1}]_{kk},\quad k = 1,\ldots,n.
\]
\end{theorem}
We omit the proof as it is identical to the proof of \cite[Theorem 2]{P95} and only uses the fact that $k$-th row and column of $(P_{k-1}\ldots P_1)L$ and $U(P_{k-1}\ldots P_1)^{-1}$ correspondingly are singletons.

\subsection{Incorporation of shift}
In order to derive the shifted version of dqds one needs to apply Gram--Schmidt process to matrices $U-\sigma L^{-1}$ and $L$:
\[
UL-\sigma I=(U-\sigma L^{-1})L=[(U-\sigma L^{-1})P^{-1}]\cdot[PL]=\widehat{L}\widehat{U}. 
\]
At first glance the additional term $-\sigma L^{-1}$ appears to spoil the derivation of $\widehat{P}$. However, since $\widehat{P} = P_{n-1}\ldots P_2P_1$ and each of $P_k^{-1}$ act only on two rows when multiplied from the right, it is not necessary to know all the entries of $L^{-1}$ in advance but only the $(k+1,k)$ entry immediately below the main diagonal. At the $k$'th step of the Gram--Schmidt process, the only change in the transformation in comparison to the unshifted case would be in the active $2\times 2$ submatrix of $(U-\sigma L^{-1})(P_1\ldots P_{k-1})^{-1}$:
\[
\begin{bmatrix}
t_k & \widehat{g}_kh_{k+1} \\
\sigma s_k & d_{k+1}-\sigma
\end{bmatrix}
\cdot
\begin{bmatrix}
1 & -\widehat{g}_kh_{k+1} \\
s_k & t_k
\end{bmatrix}
=
\begin{bmatrix}
1 & 0 \\
\widehat{s}_k & t_{k+1}
\end{bmatrix}\cdot \widehat{d}_k.
\]
This yields
\[
\begin{gathered}
t_k+s_k\widehat{g}_kh_{k+1} \equiv \widehat{d}_k,\quad\mbox{as before},\\
d_{k+1} = \widehat{s}_k\widehat{d}_k,\quad\mbox{as before},\\
-\sigma s_k\widehat{g}_kh_{k+1}+d_{k+1}t_k-\sigma t_k \equiv t_{k+1}\widehat{d}_k.
\end{gathered}
\]
The last equation simplifies to $t_{k+1} = d_{k+1}t_{k}/\widehat{d}_{k} - \sigma$, which is exactly the shifted version of dqds.

\section{dqds algorithm for companion and comrade matrices}\label{sec:companion}
In this Section we specialize dqds($\sigma$) algorithm for matrices related to polynomials, namely well-known companion and comrade matrices. In the first case the new algorithm will become a root-finding algorithm for polynomial written in monomial basis and in the second case --- in orthogonal polynomials basis.

\subsection{dqds for companion matrix}
Companion matrix has the form
\begin{equation}\label{eq:companion}
C = 
\begin{bmatrix}
-m_{n-1}  & -m_{n-2}  & \cdots & -m_1 & -m_0  \\
1     & 0     & \cdots & 0        & 0     \\
0     & 1     & \cdots & 0        & 0     \\
\cdot & \cdot & \cdot  & \cdot    & \cdot \\
0     & 0     & \cdots & 1        & 0     \\
\end{bmatrix}.
\end{equation}
and its most useful property for applications is that its characteristic polynomial $P(x)$ is encoded in the entries:
\begin{equation}\label{eq:master_poly}
P(x) = \det(xI-C) = x^n+m_{n-1}x^{n-1}+\cdots+m_1x+m_0.
\end{equation}

Let $\sigma$ be an arbitrary shift, then LU factorization of $C-\sigma I$ is
\begin{equation}\label{eq:compan_LU}
\begin{aligned}
&L = 
\begin{bmatrix}
1 \\
-\frac{H_0(\sigma)}{H_1(\sigma)} & 1 \\
& -\frac{H_1(\sigma)}{H_2(\sigma)} & 1 \\
& & \ddots & \ddots \\
& & & -\frac{H_{n-2}(\sigma)}{H_{n-1}(\sigma)} & 1
\end{bmatrix},\\
&U=
\begin{bmatrix}
-\frac{H_1(\sigma)}{H_0(\sigma)} & -\frac{m_{n-2}}{H_0(\sigma)} & \cdots & -\frac{m_1}{H_0(\sigma)} & -\frac{m_0}{H_0(\sigma)} \\
& -\frac{H_2(\sigma)}{H_1(\sigma)} & & & -\frac{m_0}{H_1(\sigma)} \\
& & \ddots & & \vdots \\
& & & -\frac{H_{n-1}(\sigma)}{H_{n-2}(\sigma)} & -\frac{m_0}{H_{n-2}(\sigma)} \\
& & & & -\frac{P(\sigma)}{H_{n-1}(\sigma)}\\
\end{bmatrix}.
\end{aligned}
\end{equation}

Where $H_k(x)$ are celebrated Horner polynomials arising in the Horner rule of evaluating $P(x)$:
\begin{equation}\label{eq:Horner}
H_0(x) = 1,\quad H_{k}(x) = x\cdot H_{k-1}(x) + m_{n-k},\quad k = 1,\ldots, n,\quad H_n(x)=P(x).
\end{equation}

Factorization \eqref{eq:compan_LU} exists for any $\sigma$ s.t. $H_k(\sigma)\neq 0$ for all $k$. In particular it exists for $\sigma=0$ if $m_k\neq 0$, $\forall k$. Quasiseparable generators $\{s_k$, $d_k$, $g_k$, $b_k$, $h_k\}$ of $L$ and $U$ factors needed for the input to the dqds Algorithm \ref{alg:dqds_Hessenberg} are given in Table \ref{tbl:shifted_compan_LU_gens}. Let us note that recursion \eqref{eq:Horner} may not be a good way of evaluating LU factorization due to the possibility of overflow, Algorithm \ref{alg:LU} provides a better way to do this.

\begin{table}[H]
\caption{Generators of $L$ and $U$ factors of a companion matrix.}\label{tbl:shifted_compan_LU_gens}
\begin{center}
\begin{tabular}[c]{>{$}c<{$}>{$}c<{$}>{$}c<{$}>{$}c<{$}>{$}c<{$}>{$}c<{$}>{$}c<{$}}
s_k & d_k & g_k & b_k & h_k\\
\midrule
-\frac{H_{k-1}(\sigma)}{H_k(\sigma)} & -\frac{H_k(\sigma)}{H_{k-1}(\sigma)} & -\frac{1}{H_{k-1}(\sigma)} & 1 & m_{n-k}\\
\bottomrule
\end{tabular}
\end{center}
\end{table}

Generators in Table \ref{tbl:shifted_compan_LU_gens} are special because $b_k=1$ for all $k$. It was shown in \cite{BEGOZ10} that quasiseparable matrix is \emph{semiseparable} if and only if there is a choice of generators such that $b_k=1$. Semiseparability means that the whole upper triangular part of matrix $U$ is an upper part of some rank-one matrix. Generators $b_k$ are preserved under iterations of Algorithm \ref{alg:dqds_Hessenberg} and, hence, semiseparable property is also preserved. We conclude that factors $L$ and $U$ stay bidiagonal and semiseparable correspondingly during the course of dqds iterations.

\subsection{dqds for comrade matrix}
Let $\{r_k(x)\}_{k=0}^{n}$ be a system of (monic) orthogonal polynomials associated with some nonnegative measure on the real line $\mathbb{R}$. As is well known, every such system satisfies a recurrence relation of the form
\[
\begin{aligned}
&r_{k+1}(x) = (x-\alpha_k)\cdot r_k(x) - \beta_k\cdot r_{k-1}(x),\quad k = 1,2,3,\ldots,\\
&r_{0}(x) = 0,\quad r_1(x) = (x-\alpha_0)\cdot r_0(x),
\end{aligned}
\]
where $\alpha_k$ and $\beta_k>0$ are unique (for a given measure and support) real constants.

Let $P(x)$ be a monic polynomial of degree $n$ represented in the basis of $\{r_k(x)\}_{k=0}^{n}$:
\[
P(x) = r_n(x) + m_{n-1}\cdot r_{n-1}(x)+\cdots+m_1\cdot r_1(x) + m_0\cdot r_0(x).
\]
Then \emph{comrade} matrix for $P(x)$ is defined as follows:
\begin{equation}\label{eq:comrade}
C = 
\begin{bmatrix}
\alpha_{n-1}-m_{n-1} & \beta_{n-1}-m_{n-2} & -m_{n-3} & \cdots & -m_0 \\
1 & \alpha_{n-2} & \beta_{n-2} & &  \\
& 1 & \ddots & \ddots & \\
& & \ddots & \alpha_1 & \beta_1 \\
& & & 1 & \alpha_0 \\
\end{bmatrix}.
\end{equation}

For an arbitrary value of the shift $\sigma$, the LU factorization of shifted comrade matrix $C-\sigma I$ is 
\begin{equation}\label{eq:comrade_LU}
\begin{aligned}
&L = 
\begin{bmatrix}
1 \\
-\frac{H_0(\sigma)}{H_1(\sigma)} & 1 \\
& -\frac{H_1(\sigma)}{H_2(\sigma)} & 1 \\
& & \ddots & \ddots \\
& & & -\frac{H_{n-2}(\sigma)}{H_{n-1}(\sigma)} & 1
\end{bmatrix},\\
&U=
\begin{bmatrix}
-\frac{H_1(\sigma)}{H_0(\sigma)} & \beta_{n-1}-\frac{m_{n-2}}{H_0(\sigma)} & \cdots & -\frac{m_1}{H_0(\sigma)} & -\frac{m_0}{H_0(\sigma)} \\
& -\frac{H_2(\sigma)}{H_1(\sigma)} & \ddots & & -\frac{m_0}{H_1(\sigma)} \\
& & \ddots & \ddots & \vdots \\
& & & -\frac{H_{n-1}(\sigma)}{H_{n-2}(\sigma)} & \beta_1-\frac{m_0}{H_{n-2}(\sigma)} \\
& & & & -\frac{P(\sigma)}{H_{n-1}(\sigma)}\\
\end{bmatrix},
\end{aligned}
\end{equation}
where $H_k(x)$ are generalized Horner polynomials defined via the recurrence relation (Clenshaw algorithm \cite{C55}):
\begin{equation}\label{eq:genHorner}
\begin{aligned}
&H_{k}(x) = (x-\alpha_{n-k})\cdot H_{k-1}(x) - \beta_{n-k+1}H_{k-2}(x) + m_{n-k},\quad k = 2,\ldots, n,\\
&H_0(x) = 1,\quad H_1(x) = x - \alpha_{n-1} + m_{n-1},\quad H_n(x)=P(x).
\end{aligned}
\end{equation}

LU factorization \eqref{eq:comrade_LU} exists if and only if $H_k(\sigma)\neq 0$ for all $k=1,\ldots,n$. Matrix $U$ is quasiseparable of order $2$ and, therefore, its generators $g_k$, $b_k$, $h_k$ are matrices of appropriate sizes. All generators of $L$ and $U$ factors are given in Table \ref{tbl:shifted_comrade_LU_gens}.

\begin{table}[H]
\begin{center}
\caption{Generators of $L$ and $U$ factors of a shifted comrade matrix.}\label{tbl:shifted_comrade_LU_gens}
\begin{tabular}[c]{>{$}c<{$}>{$}c<{$}>{$}c<{$}>{$}c<{$}>{$}c<{$}>{$}c<{$}}
s_k & d_k & g_k & b_k & h_k\\
\midrule
-\frac{H_{k-1}(\sigma)}{H_k(\sigma)} &
-\frac{H_k(\sigma)}{H_{k-1}(\sigma)} &
\begin{bmatrix}1&-\frac{1}{H_{k-1}(\sigma)}\end{bmatrix} &
\begin{bmatrix}0&0\\0&1\end{bmatrix} & 
\begin{bmatrix}\beta_{n-k+1}&m_{n-k}\end{bmatrix}^T \\
\bottomrule
\end{tabular}
\end{center}
\end{table}

Generators in Table \ref{tbl:shifted_comrade_LU_gens} are input parameters to the dqds($\sigma$) Algorithm \ref{alg:dqds_Hessenberg} that one can run with suitable shifts to compute roots of polynomial $P(x)$. Due to the special form of generators of $U$ and some conservation properties of Algorithm \ref{alg:dqds_Hessenberg}, matrix $U$ stays the sum of a bidiagonal and semiseparable of order one at every step of dqds iterations.

\subsection{dqds for fellow matrix}
One way to evaluate roots of a polynomial represented in the basis of polynomials orthogonal on the unit circle (also called Szeg\"o polynomials), is to form a so-called \emph{fellow} matrix. Fellow matrix is also Hessenberg and is an analog of comrade matrix. Hovewer, in this case its upper triangular part is dense. It is known that fellow matrix is quasiseparable of order $2$ (see, for instance, \cite{BEGOTZ09}) and, therefore, Algorithm \ref{alg:dqds_Hessenberg} applies to it. We do not derive initial LU factorization of a fellow matrix here but, if needed, it can be done by analogy with companion and comrade cases.

\section{Numerical results}\label{sec:results}
In this section we present results of numerical tests with the new dqds Algorithm \ref{alg:dqds_Hessenberg} implemented in MATLAB (2010b, The MathWorks Inc., Natick, MA). We compare performances of the following algorithms for finding roots of polynomials:
\begin{itemize}
\item \textsf{qs\_dqds} --- qusiseparable dqds Algorithm \ref{alg:dqds_Hessenberg};
\item \textsf{compan\_qr} --- companion QR algorithm described in \cite{BBEGG10}. The implementation of this algorithm was generously provided by Paola Boito and is also accessible at \url{www.unilim.fr/pages_perso/paola.boito/software.html};
\item \textsf{eig} --- MATLAB's main eigenvalue routine. When applied to a comrade matrix \eqref{eq:comrade} can be used to find roots of polynomials represented in orthogonal polynomials' bases.
\item \textsf{roots} --- MATLAB's main routine for finding polynomial's roots. Same as \textsf{eig} applied to a companion matrix.
\end{itemize}

Also in the tables of this section we will use the following notations:
\begin{itemize}
\item $n$ --- polynomial's degree;
\item ni --- number of iterations per root used by an algorithm;
\item $\varepsilon$ --- relative accuracy of the computed roots evaluated as $\varepsilon = \max_i\frac{|x_i-\widehat{x}_i|}{|x_i|}$, where $x_i$ and $\widehat{x}_i$ are exact and computed roots, respectively.
\end{itemize}

It is a common strategy to perform diagonal scaling (balancing) of a matrix before applying an eigenvalue algorithm (see, for instance, \cite{PR69}). Balancing can greatly improve accuracy of the computed eigenvalues. Unfortunately, algorithm \textsf{compan\_qr} does not allow balancing because it represents initial companion matrix and all the subsequent iterates as unitary plus rank-one matrices and balancing destroys this structure. For completeness we present results of all the remaining algorithms with and without balancing. For balancing algorithm we simply used MATLAB's \textsf{balance} routine.

Coefficients of test polynomials in monomials basis were computed in MATLAB from their roots using the $\textsf{poly}$ command. Coefficients in orthogonal polynomials bases were computed from the corresponding coefficients in monomials basis (details of this algorithm can be found, for instance, in \cite{BEGOTZ09}).

Details of our implementation of \textsf{qs\_dqds} algorithm are as follows. The initial LU factorization was computed with zero shift. Value of the shift $\sigma$ at every consecutive iteration was chosen to be the current iterate's $(k,k)$ entry (when computing $k$-th eigenvalue). We used a simple deflation criterion $|A_{k,k-1}|<\varepsilon |A_{k,k}|$ with $\varepsilon=10^{-12}$. Since, at the moment \textsf{qs\_dqds} algorithm is only developed for single shifts, in our experiments we have considered polynomials with real roots. Current implementation and the script that reproduces results of this section are avalable at \url{www.maths.ed.ac.uk/~pzhlobic/software.shtml}.

\textbf{Test 1:} Wilkinson's first polynomial $p(x)=\prod\limits_{i=1}^n(x-i)$.\\
It is known due J.H.~Wilkinson \cite{W59} that large roots of this polynomial are extremely sensitive to the perturbation of coefficients. Hence, computing roots of this polynomial is a hard test for any root-finding algorithm and we can't expect any such algorithm to compute these roots accurately for a polynomial of high degree. Results of numerical tests presented in Table \ref{tbl:wp1} demonstrate that algorithm \textsf{qs\_dqds} always delivers same or higher relative accuracy and the number of iterations per root is quite moderate (3.3 on average). In the companion matrix case our implementation modifies parameters ``in place'' and, therefore, uses only $4n$ storage and $10n$ arithmetic operations per iteration.

\begin{table}[!ht]
\begin{center}
\caption{Wilkinson's first polynomial.}\label{tbl:wp1}
\begin{tabular}[c]{r|rr|r|rr|rr|r}& \multicolumn{3}{c|}{with balancing} & \multicolumn{5}{c}{without balancing}\\\midrule& \multicolumn{2}{c|}{\textsf{qs\_dqds}} & \multicolumn{1}{c|}{\textsf{roots}} & \multicolumn{2}{c|}{\textsf{qs\_dqds}} & \multicolumn{2}{c|}{\textsf{\textsf{compan\_qr}}} & \multicolumn{1}{c}{\textsf{roots}}\\$n$ & \multicolumn{1}{c}{$\varepsilon$} & \multicolumn{1}{c|}{ni} & \multicolumn{1}{c|}{$\varepsilon$} & \multicolumn{1}{c}{$\varepsilon$} & \multicolumn{1}{c|}{ni} & \multicolumn{1}{c}{$\varepsilon$} & \multicolumn{1}{c|}{ni} & \multicolumn{1}{c}{$\varepsilon$}\\\midrule
10 & 1.6e-11 & 3.4 & 5.6e-10 & 2.1e-11 & 3.3 & 7.1e-11 & 4.0 & 5.0e-11 \\
11 & 9.3e-11 & 3.4 & 2.4e-09 & 7.5e-11 & 3.3 & 1.3e-10 & 4.1 & 1.1e-09 \\
12 & 2.5e-09 & 3.4 & 2.7e-08 & 2.4e-09 & 3.3 & 1.6e-09 & 4.1 & 2.9e-09 \\
13 & 1.2e-08 & 3.4 & 4.6e-08 & 1.2e-08 & 3.3 & 2.9e-09 & 4.0 & 3.6e-08 \\
14 & 1.1e-08 & 3.4 & 3.8e-07 & 1.1e-08 & 3.3 & 3.7e-08 & 4.0 & 5.6e-07 \\
15 & 7.3e-08 & 3.4 & 2.5e-06 & 7.3e-08 & 3.3 & 2.3e-07 & 3.8 & 1.0e-06 \\
16 & 8.8e-08 & 3.4 & 8.8e-06 & 8.8e-08 & 3.3 & 5.2e-07 & 4.0 & 1.4e-06 \\
17 & 7.6e-06 & 3.5 & 4.5e-05 & 7.6e-06 & 3.4 & 4.1e-06 & 4.0 & 3.2e-05 \\
18 & 2.2e-05 & 3.4 & 2.7e-04 & 2.2e-05 & 3.3 & 2.6e-05 & 3.8 & 2.4e-04 \\
19 & 1.2e-04 & 3.4 & 6.1e-04 & 1.2e-04 & 3.3 & 1.1e-04 & 3.8 & 2.8e-04 \\
20 & 9.4e-04 & 3.4 & 4.7e-03 & 9.4e-04 & 3.4 & 4.7e-04 & 3.9 & 4.5e-03 \\
\bottomrule
\end{tabular}
\end{center}
\end{table}

\textbf{Test 2:} Reversed Wilkinson's first polynomial $p(x)=\prod\limits_{i=1}^n(x-1/i)$.\\
This polynomial happened to be a harder test for the algorithms. Roots computed by \textsf{qs\_dqds} were from $10$ to $10^4$ times more accurate than those computed by MATLAB's \textsf{roots} function. Also note that surprisingly balancing did not have any influence on the accuracy of roots computed by \textsf{qs\_dqds}. Structured QR iterations algorithm \textsf{compan\_qr} was the least accurate on this test.

\begin{table}[!ht]
\begin{center}
\caption{Reversed Wilkinson's first polynomial.}\label{tbl:rwp1}
\begin{tabular}[c]{r|rr|r|rr|rr|r}& \multicolumn{3}{c|}{with balancing} & \multicolumn{5}{c}{without balancing}\\\midrule& \multicolumn{2}{c|}{\textsf{qs\_dqds}} & \multicolumn{1}{c|}{\textsf{roots}} & \multicolumn{2}{c|}{\textsf{qs\_dqds}} & \multicolumn{2}{c|}{\textsf{\textsf{compan\_qr}}} & \multicolumn{1}{c}{\textsf{roots}}\\$n$ & \multicolumn{1}{c}{$\varepsilon$} & \multicolumn{1}{c|}{ni} & \multicolumn{1}{c|}{$\varepsilon$} & \multicolumn{1}{c}{$\varepsilon$} & \multicolumn{1}{c|}{ni} & \multicolumn{1}{c}{$\varepsilon$} & \multicolumn{1}{c|}{ni} & \multicolumn{1}{c}{$\varepsilon$}\\\midrule
10 & 1.6e-10 & 3.3 & 2.4e-09 & 1.6e-10 & 3.3 & 3.0e-07 & 4.8 & 6.7e-07 \\
11 & 1.1e-10 & 3.3 & 1.8e-08 & 1.1e-10 & 3.3 & 1.5e-04 & 4.2 & 1.5e-04 \\
12 & 1.6e-09 & 3.3 & 1.8e-07 & 1.6e-09 & 3.3 & 2.7e-02 & 4.8 & 1.3e-04 \\
13 & 3.0e-09 & 3.2 & 1.4e-06 & 3.0e-09 & 3.2 & 1.1e-01 & 7.8 & 1.4e-01 \\
14 & 5.5e-08 & 3.3 & 1.0e-05 & 5.5e-08 & 3.3 & 1.3e-01 & 5.7 & 2.4e-01 \\
15 & 6.6e-08 & 3.3 & 1.9e-04 & 6.6e-08 & 3.3 & 4.6e-01 & 7.0 & 5.1e-01 \\
16 & 2.0e-06 & 3.2 & 1.0e-03 & 2.0e-06 & 3.3 & 4.8e-01 & 6.9 & 6.8e-01 \\
17 & 6.5e-06 & 3.2 & 1.3e-02 & 6.5e-06 & 3.2 & 1.1e+00 & 7.8 & 7.3e-01 \\
18 & 5.3e-05 & 3.2 & 1.6e-01 & 5.3e-05 & 3.3 & 1.1e+00 & 7.1 & 1.1e+00 \\
19 & 1.5e-04 & 3.2 & 1.8e-01 & 1.5e-04 & 3.2 & 1.4e+00 & 10.3 & 1.4e+00 \\
20 & 3.7e-03 & 3.2 & 2.7e-01 & 3.7e-03 & 3.3 & 1.3e+00 & 7.2 & 1.9e+00 \\
\bottomrule
\end{tabular}
\end{center}
\end{table}

\textbf{Test 3:} Wilkinson's second polynomial $p(x)=\prod\limits_{i=1}^n(x-0.6^i)$. \\
In contrast to the Wilkinson's first polynomial case, roots of this polynomial are stable with respect to small perturbations of coefficients, although the ratio of the largest to the smallest root grows exponentially. Results of numerical tests with this polynomial are presented in Table \ref{tbl:wp2}. Both \textsf{roots} and \textsf{compan\_qr} algorithms failed to compute roots accurately for polynomials of degrees higher than 30. In the case of \textsf{roots} even balancing did not improve the situation. On the contrary \textsf{qs\_dqds} maintained high accuracy of the computed roots for polynomials of degrees up to 50.

\begin{table}[!ht]
\begin{center}
\caption{Wilkinson's second polynomial.}\label{tbl:wp2}
\begin{tabular}[c]{r|rr|r|rr|rr|r}& \multicolumn{3}{c|}{with balancing} & \multicolumn{5}{c}{without balancing}\\\midrule& \multicolumn{2}{c|}{\textsf{qs\_dqds}} & \multicolumn{1}{c|}{\textsf{roots}} & \multicolumn{2}{c|}{\textsf{qs\_dqds}} & \multicolumn{2}{c|}{\textsf{\textsf{compan\_qr}}} & \multicolumn{1}{c}{\textsf{roots}}\\$n$ & \multicolumn{1}{c}{$\varepsilon$} & \multicolumn{1}{c|}{ni} & \multicolumn{1}{c|}{$\varepsilon$} & \multicolumn{1}{c}{$\varepsilon$} & \multicolumn{1}{c|}{ni} & \multicolumn{1}{c}{$\varepsilon$} & \multicolumn{1}{c|}{ni} & \multicolumn{1}{c}{$\varepsilon$}\\\midrule
10 & 4.6e-13 & 2.7 & 3.7e-13 & 4.8e-14 & 2.8 & 1.8e-03 & 3.3 & 1.7e-03 \\
20 & 8.2e-13 & 2.4 & 2.8e-09 & 6.4e-14 & 2.5 & 2.1e+00 & 7.8 & 3.9e+00 \\
30 & 8.1e-13 & 2.0 & 3.2e-04 & 2.1e-13 & 2.2 & 1.3e+04 & 11.4 & 6.9e+05 \\
40 & 8.4e-13 & 1.8 & 2.4e+00 & 1.8e-13 & 2.0 & 2.1e+07 & 7.7 & 1.7e+08 \\
50 & 1.6e-11 & 1.6 & 1.3e+01 & 2.5e-13 & 1.9 & 7.6e+09 & 7.6 & 3.4e+10 \\
\bottomrule
\end{tabular}
\end{center}
\end{table}

\textbf{Test 4:} To demonstrate the robustness of the proposed algorithm we conducted 1000 random tests with polynomials having roots randomly distributed as $X \cdot 10^{5Y}$, where $X$ and $Y$ are random variables uniformly distributed on $[-1,1]$. Such polynomials have both positive and negative roots of very different magnitudes. We report mean and standard deviation of the relative error in the tests in Table \ref{tbl:randp}. Tests indicate that \textsf{qs\_dqds} becomes more accurate (on average) than \textsf{roots} algorithm as polynomial's degree increases.

\begin{table}[!ht]
\begin{center}
\caption{Polynomial with roots distributed as $X \cdot 10^{5Y}$, $X,Y\sim U[-1,1]$.}\label{tbl:randp}
\begin{tabular}[c]{r|rr|r|rr|r}& \multicolumn{3}{c|}{with balancing} & \multicolumn{3}{c}{without balancing}\\\midrule& \multicolumn{2}{c|}{\textsf{qs\_dqds}} & \multicolumn{1}{c|}{\textsf{roots}} & \multicolumn{2}{c|}{\textsf{qs\_dqds}} & \multicolumn{1}{c}{\textsf{roots}}\\$n$ & \multicolumn{1}{c}{$\varepsilon$} & \multicolumn{1}{c|}{ni} & \multicolumn{1}{c|}{$\varepsilon$} & \multicolumn{1}{c}{$\varepsilon$} & \multicolumn{1}{c|}{ni} & \multicolumn{1}{c}{$\varepsilon$}\\\midrule
10 & e-13$\pm$e-12 & 2.2$\pm$0.3 & e-12$\pm$e-11 & e-09$\pm$e-08 & 2.3$\pm$0.3 & e+07$\pm$e+08 \\
20 & e-12$\pm$e-11 & 2.4$\pm$0.2 & e-09$\pm$e-08 & e-08$\pm$e-07 & 2.4$\pm$0.2 & e+07$\pm$e+08 \\
30 & e-11$\pm$e-10 & 2.4$\pm$0.2 & e-07$\pm$e-06 & e-08$\pm$e-07 & 2.4$\pm$0.2 & e+06$\pm$e+07 \\
\bottomrule
\end{tabular}
\end{center}
\end{table}

\textbf{Test 5:} Algorithm \textsf{qs\_dqds} can be applied to polynomials represented in a basis of arbitrary orthogonal polynomials. We confirm this by presenting results of numerical tests with Wilkinson's polynomial represented in the basis of Chebyshev polynomials of the second kind. Chebyshev polynomials were scaled to the interval $[a-1,b+1]$, where $a$ and $b$ is the smallest and the largest root of the polynomial of interest, respectively. We compared performance of the new algorithm with MATLAB's \textsf{eig} routine. Results are presented in Table \ref{tbl:cwp1}.

\begin{table}[!ht]
\begin{center}
\caption{Wilkinson's first polynomial in Chebyshev basis.}\label{tbl:cwp1}
\begin{tabular}[c]{r|rr|r|rr|r}& \multicolumn{3}{c|}{with balancing} & \multicolumn{3}{c}{without balancing}\\\midrule& \multicolumn{2}{c|}{\textsf{qs\_dqds}} & \multicolumn{1}{c|}{\textsf{eig}} & \multicolumn{2}{c|}{\textsf{qs\_dqds}} & \multicolumn{1}{c}{\textsf{eig}}\\$n$ & \multicolumn{1}{c}{$\varepsilon$} & \multicolumn{1}{c|}{ni} & \multicolumn{1}{c|}{$\varepsilon$} & \multicolumn{1}{c}{$\varepsilon$} & \multicolumn{1}{c|}{ni} & \multicolumn{1}{c}{$\varepsilon$}\\\midrule
10 & 6.5e-13 & 4.1 & 4.4e-15 & 6.5e-13 & 4.1 & 1.9e-15 \\
12 & 1.5e-12 & 3.9 & 6.5e-15 & 1.5e-12 & 3.9 & 3.4e-15 \\
14 & 5.9e-13 & 3.9 & 1.6e-13 & 5.9e-13 & 3.9 & 1.5e-13 \\
16 & 7.0e-12 & 4.3 & 6.9e-12 & 7.0e-12 & 4.3 & 6.9e-12 \\
18 & 5.8e-11 & 4.4 & 5.8e-11 & 5.8e-11 & 4.4 & 5.8e-11 \\
20 & 1.8e-09 & 4.8 & 1.8e-09 & 1.8e-09 & 4.8 & 1.8e-09 \\
22 & 7.0e-09 & 4.7 & 7.0e-09 & 7.0e-09 & 5.1 & 7.0e-09 \\
24 & 1.0e-07 & 5.0 & 1.0e-07 & 1.0e-07 & 5.0 & 1.0e-07 \\
26 & 7.8e-07 & 4.1 & 7.8e-07 & 7.8e-07 & 4.1 & 7.8e-07 \\
28 & 1.1e-05 & 3.5 & 1.1e-05 & 1.1e-05 & 3.5 & 1.1e-05 \\
30 & 3.5e-05 & 4.3 & 3.5e-05 & 3.5e-05 & 4.3 & 3.5e-05 \\
\bottomrule
\end{tabular}
\end{center}
\end{table}

\section*{Conclusion}
We extended dqds algorithm of Fernando and Parlett \cite{FP94,P95} from tridiagonal to quasiseparable matrices only in the Hessenberg case but qds algorithm to all quasiseparable matrices. New algorithms use just $\mathcal{O}(n)$ operations per iteration. Hessenberg quasiseparable matrices include, among others, important companion and comrade matrices and, hence, can be used to find roots of polynomials represented in different polynomial bases, e.g. monomials, orthogonal polynomials. Our findings are supported by preliminary numerical experiments with polynomials having real roots. In our implementation we used naive shift strategy and single shifts only. There are numerous extensions and improvements of the algorithm left for future research including the development of more sophisticated shift strategies for better robustness and stability as well as double shifts for handling complex roots.

\section*{Acknowledgement}
The author is grateful to the Editor and two anonymous referees for very constructive comments and suggestions that helped to improve the paper considerably. Special thanks are also extended to Paola Boito for providing implementation of the \textsf{compan\_qr} algorithm as well as the useful remarks on balancing of companion matrix.

\bibliographystyle{siam}

\begin{thebibliography}{10}

\bibitem{BEGOZ10}
{\sc T.~Bella, Y.~Eidelman, I.~Gohberg, V.~Olshevsky, and P.~Zhlobich}, {\em
  Classifications of recurrence relations via subclasses of
  ({H},m)-quasiseparable matrices}, in Numerical Linear Algebra in Signals,
  Systems and Control, vol.~XV of Lecture Notes in Electrical Engineering,
  {Springer--Verlag}, 2011, pp.~23--53.

\bibitem{BEGOTZ09}
{\sc T.~Bella, V.~Olshevsky, P.~Zhlobich, Y.~Eidelman, I.~Gohberg, and
  E.~Tyrtyshnikov}, {\em {A Traub-like algorithm for
  Hessenberg-quasiseparable-Vandermonde matrices of arbitrary order}}, in
  Numerical Methods for Structured Matrices and Applications, D.~Bini,
  V.~Mehrmann, V.~Olshevsky, E.~Tyrtyshnikov, and M.~{Van Barel}, eds.,
  vol.~199 of Operator Theory: Advances and Applications, Birkh{\"a}user Basel,
  2010, pp.~127--154.

\bibitem{BBD11}
{\sc R.~Bevilacqua, E.~Bozzo, and G.~Del~Corso}, {\em {qd-type methods for
  quasiseparable matrices}}, SIAM Journal on Matrix Analysis and Applications,
  32 (2011), pp.~722--747.

\bibitem{BBEGG10}
{\sc D.~Bini, P.~Boito, Y.~Eidelman, L.~Gemignani, and I.~Gohberg}, {\em {A
  fast implicit QR eigenvalue algorithm for companion matrices}}, Linear
  Algebra and its Applications, 432 (2010), pp.~2006--2031.

\bibitem{BDG04}
{\sc D.~Bini, F.~Daddi, and L.~Gemignani}, {\em On the shifted {QR} iteration
  applied to companion matrices}, Electronic Transactions on Numerical
  Analysis, 18 (2004), pp.~137--152.

\bibitem{BEGG07}
{\sc D.~Bini, Y.~Eidelman, L.~Gemignani, and I.~Gohberg}, {\em {Fast QR
  eigenvalue algorithms for Hessenberg matrices which are rank-one
  perturbations of unitary matrices}}, SIAM Journal on Matrix Analysis and
  Applications, 29 (2007), pp.~566--585.

\bibitem{BGP05}
{\sc D.~Bini, L.~Gemignani, and V.~Pan}, {\em {Fast and stable QR eigenvalue
  algorithms for generalized companion matrices and secular equations}},
  Numerische Mathematik, 100 (2005), pp.~373--408.

\bibitem{CGXZ08}
{\sc S.~Chandrasekaran, M.~Gu, J.~Xia, and J.~Zhu}, {\em {A fast QR algorithm
  for companion matrices}}, Recent Advances in Matrix and Operator Theory, 179
  (2008), pp.~111--143.

\bibitem{C55}
{\sc C.~Clenshaw}, {\em {A note on the summation of Chebyshev series}}, MTAC,
  v. 9, 17 (1955), pp.~118--120.

\bibitem{DV07}
{\sc S.~Delvaux and M.~Van~Barel}, {\em A {G}ivens-weight representation for
  rank structured matrices}, SIAM Journal on Matrix Analysis and Applications,
  29 (2007), pp.~1147--1170.

\bibitem{DK90}
{\sc J.~Demmel and  W.~Kahan}, {\em Accurate singular values of bidiagonal matrices}, SIAM Journal on Scientific and Statistical Computing, 11 (1990), pp.~873--912.

\bibitem{DOZ11}
{\sc F.~Dopico, V.~Olshevsky, and P.~Zhlobich}, {\em Stability of qr-based
  system solvers for a subclass of quasiseparable rank one matrices}, Mathematics of Computation, in press
   (2011).

\bibitem{EG99a}
{\sc Y.~Eidelman and I.~Gohberg}, {\em On a new class of structured matrices},
  Integral Equations and Operator Theory, 34 (1999), pp.~293--324.

\bibitem{EG05}
\leavevmode\vrule height 2pt depth -1.6pt width 23pt, {\em On generators of
  quasiseparable finite block matrices}, Calcolo, 42 (2005), pp.~187--214.

\bibitem{EGO05b}
{\sc Y.~Eidelman, I.~Gohberg, and V.~Olshevsky}, {\em {The QR iteration method
  for Hermitian quasiseparable matrices of an arbitrary order}}, Linear Algebra
  and its Applications, 404 (2005), pp.~305--324.

\bibitem{FP94}
{\sc K.~Fernando and B.~Parlett}, {\em Accurate singular values and
  differential qd algorithms}, Numerische Mathematik, 67 (1994), pp.~191--229.

\bibitem{GKK85}
{\sc I.~Gohberg, T.~Kailath, and I.~Koltracht}, {\em Linear complexity
  algorithms for semiseparable matrices}, Integral Equations and Operator
  Theory, 8 (1985), pp.~780--804.

\bibitem{GV93}
{\sc G.~H. Golub and C.~F. {Van Loan}}, {\em {Matrix Computations}}, Johns
  Hopkins Series in the Mathematical Sciences, Johns Hopkins University Press,
  Baltimore and London, second~ed., 1993.

\bibitem{PR69}
{\sc B.~Parlett and C.~Reinsch}, {\em Balancing a matrix for calculation of eigenvalues and eigenvectors}, Numerische Mathematik, 13 (1969), pp.~292--304.

\bibitem{P95}
{\sc B.~Parlett}, {\em The new qd algorithms}, Acta Numerica, 4 (1995),
  pp.~459--491.

\bibitem{PVV08}
{\sc B.~Plestenjak, M.~Van~Barel, and E.~Van~Camp}, {\em {A Cholesky LR
  algorithm for the positive definite symmetric diagonal-plus-semiseparable
  eigenproblem}}, Linear Algebra and its Applications, 428 (2008),
  pp.~586--599.

\bibitem{VVM08}
{\sc R.~Vandebril, M.~Van~Barel, and N.~Mastronardi}, {\em Matrix computations and semiseparable matrices}, Vol.~I: Linear Systems, Johns Hopkins University Press, 2008.

\bibitem{VVM05b}
\leavevmode\vrule height 2pt depth -1.6pt width 23pt, {\em {An implicit QR
  algorithm for symmetric semiseparable matrices}}, Numerical Linear Algebra
  with Applications, 12 (2005), pp.~625--658.

\bibitem{W59}
{\sc J.H.~Wilkinson}, {\em The evaluation of the zeros of ill-conditioned polynomials. Part I}, Numerische Mathematik, 1(1959) pp.~150--166.

\bibitem{XCGL09}
{\sc J.~Xia, S.~Chandrasekaran, M.~Gu, and X.~Li}, {\em Superfast multifrontal
  method for large structured linear systems of equations}, SIAM Journal on
  Matrix Analysis and Applications, 31 (2009), pp.~1382--1411.

\end{thebibliography}

\end{document}